\documentclass[11pt]{article}
\usepackage{amsthm,amssymb,amsmath,amscd}
\usepackage[all]{xy}
\usepackage{mathrsfs}
\usepackage{multirow}

\usepackage{mathptm}

\newtheorem{Def}{Definition}[section]
\newtheorem{Prop}[Def]{Proposition}
\newtheorem{Theo}[Def]{Theorem}
\newtheorem{Lem}[Def]{Lemma}
\newtheorem{Koro}[Def]{Corollary}
\newtheorem{Bsp}[Def]{Example}

\newcommand{\Tor}{{\rm Tor}}

\newcommand{\add}{{\rm add}}

\newcommand{\Hom}{{\rm Hom }}
\newcommand{\pd}{{\rm pd }}

\newcommand{\rad}{{\rm rad}}

\newcommand{\End}{{\rm End}}
\newcommand{\Ext}{{\rm Ext}}

\newcommand{\modcat}[1]{#1\mbox{{\rm -mod}}}

\newcommand{\opp}{^{\rm op}}
\newcommand{\otimesL}{\otimes^{\rm\bf L}}

\newcommand{\rHom}{{\rm\bf R}{\rm Hom}}

\newcommand{\gd}{{\rm gl.dim}}

\newcommand{\lra}{\longrightarrow}
\newcommand{\ra}{\rightarrow}

\begin{document}

{\Large \bf
\begin{center}
Ringel's contributions to quasi-hereditary algebras\\
\end{center}}
\centerline{Dedicated to Claus Michael Ringel's  80th Birthday}

\bigskip
\centerline{{\bf Changchang Xi}}

\begin{abstract}
Quasi-hereditary algebras were introduced by Cline, Parshall and Scott to describe the highest weight categories of representations of semisimple Lie algebras and algebraic groups by the module categories of finite-dimensional algebras. Since then a lot of homological, structural and categorical properties of quasi-hereditary algebras have been discovered. This class of algebras seems quite common and occurs in many branches of mathematics. There are lots of important works on the subject. In this note we mainly survey some of Claus Michael Ringel's works or works jointly with his collaborators on quasi-hereditary algebras. Also, some of related works and recent developments on quasi-hereditary algebras are mentioned.
\end{abstract}

{\footnotesize\tableofcontents\label{contents}}

\renewcommand{\thefootnote}{\alph{footnote}}
\setcounter{footnote}{-1} 
\renewcommand{\thefootnote}{\alph{footnote}}
\setcounter{footnote}{-1} 

\section{Introduction}
Quasi-hereditary algebras were introduced by Cline, Parshall and Scott in order
to study highest weight categories in the representation theory of semisimple
complex Lie algebras \cite{bgg} and  semisimple algebraic groups (see \cite{s,cps1988b, ps1988}). They axiomatize highest weight categories and show that the highest weight categories with a finite number of weights are in fact just the module categories of finite-dimensional quasi-hereditary algebras. Since then a lot of classes of algebras have been proved to be quasi-hereditary in many mathematical branches including representation theory of algebras, algebraic groups, category $\mathcal{O}$ in Lie theory, quantum groups and statistic mechanics, and there has been a large variety of contributions to the beautiful theory of quasi-hereditary algebras in different aspects. Especially, a lot of significant works have been done by Claus Michael Ringel or jointly with Vlastimil Dlab in ring theory and representation theory of algebras.

In this note, I will try to survey some strick results and creative ideas of Ringel (or jointly with his collaborators) on the study of quasi-hereditary algebras and related works.

The contents of this note are outlined as follows. In Section \ref{sect2} we give the ring theoretical definition of quasi-hereditary rings (algebras) and examples. We summary basic properties of quasi-hereditary algebras, which were developed by Dlab and Ringel. At the end of this section we mention some generalizations of the notion of quasi-hereditary algebras, and concentrate us mainly on standardly stratified algebras. In Section \ref{sect3} we present an influent theory of Ringel's and Dlab--Ringel's module theoretical approach to quasi-hereditary algebras. We focuses on the study of the category of good modules. Some generalizations such as standardization systems are included. In Section \ref{sect3+} we survey Koenig's decomposition theory of quasi-hereditary algebras as exact Borel subalgebras and Delta-subalgebras. Some recent results on the category of good modules are also mentioned. In Section \ref{sect4} we introduce the systematic construction of quasi-hereditary algebras by Dlab and Ringel, which theoretically produces all quasi-hereditary rings (algebras). In this section we also quote a construction of quasi-hereditary algebras from self-injective algebras. This construction is due to Dlab, Heath and Marko. Such a class of algebras may potentially have applications in understanding homological conjectures in the representation theory and homological algebra of finite-dimensional algebras.

Because of my limited knowledge and understanding, this survey is unable to exhaust all of the interesting theories and results on quasi-hereditary algebras in the literature, as expected usually. I am sorry for my choices of the contents in this note, but all of the reader's comments and suggestive critics are welcome.

\section{Ring theoretical definition of quasi-hereditary algebras\label{sect2}}

There are two ways to define quasi-hereditary rings (or algebras). One is in terms of ring theoretical language as it was done historically. The other is in terms of  module theoretical language.

In this section we give the first one and survey basic properties of quasi-hereditary algebras and heredity ideals.

By algebras we mean Artin algebras or finite-dimensional algebras, and by modules we mean left modules.

For a semiprimary ring (or Artin algebra) $A$, we denote by $A\modcat$ the category of finitely generated $A$-modules. The composition of two homomorphisms $f:X\to Y$ and $g: Y\to Z$ of $A$-modules is denoted by $fg: X\to Z$. Our notation and terminology are standard, for non-explained ones, we refer to \cite{ringel1099}.

Recall that a ring $A$ with identity is called a \emph{semiprimary }ring if the Jacobson radical $N$ of $A$ is nilpotent and $A/N$  is a semisimple artinian ring. An ideal $I$ of $A$ is idempotent if $I^2=I$. Note that an ideal $I$ of a semiprimary ring $A$ is idempotent if and only if $I=AeA$ for an idempotent $e^2=e\in A$ (see \cite[Statement 6]{dr1989}).

\subsection{Definitions and examples\label{def-exa}}

Let $A$ be an associative ring with identity, and $N$ its Jacobson radical.

The notion of heredity ideals play an important role in defining quasi-hereditary algebras (rings).

\begin{Def}\label{def2.1}{\rm \cite{cps1988b}} $(1)$ An ideal $I$ of $A$ is called a
\emph{heredity ideal} if $I^2=I, INI=0$ and $_AI$ is a projective $A$-module.

$(2)$ A semiprimary ring $A$ is said to be \emph{quasi-hereditary} if
there is a chain of idempotent ideals
$$ (*) \quad 0=I_0\subseteq I_1\subseteq \cdots\subseteq I_n=A$$
of $A$ such that $I_j/I_{j-1}$ is a heredity ideal of
$A/I_{j-1}$ for all $1\le j\le n$. Such a chain of idempotent ideals is called a \emph{heredity chain} of $A$, denoted by $\mathscr{I}=(I_i)_{0\le i\le n}$.
\end{Def}

A semiprimary ring $A$ is quasi-hereditary if and only if the opposite ring $A^{\opp}$ of $A$ is quasi-hereditary (see \cite{ps1988} and \cite{dr1989a}).

Given  a quasi-hereditary ring $A$ with a heredity chain $\mathscr{I}=(I_i)_{0\le i\le n}$ and an $A$-module $_AM$, there is a canonical filtration of submodules of $M$:
$$ 0=I_0M\subseteq I_1M\subseteq \cdots\subseteq I_n=AM=M.$$
This filtration is called an \emph{$\mathscr{I}$-filtration} of $M$. It is called a \emph{good} $\mathscr{I}$-filtration if $I_iM/I_{i-1}M$ is a projective $A/I_{i-1}$-module for all $1\le i\le n$.

From Definition \ref{def2.1} we see immediately that the projectivity of an idempotent ideal $AeA$ of $A$ is of great importance for understanding heredity ideals. To know if an idempotent ideal $_AI$ of $A$ is projective, Dlab--Ringel provides us with the following criterion.

\begin{Lem} {\rm \cite{dr1989}} Let $e$ be an idempotent of a ring $R$. If the left module
$_R(ReR)$ or the right module $(ReR)_R$ is projective, then the multiplication map
$\mu: Re\otimes_{eRe}eR\ra ReR, \; xe\otimes ey\mapsto xey,$
is bijective. Conversely, assume that $R$ is semiprimary with the Jacobson radical $N$ such that
$eNe=0$. Then both modules $(ReR)_R$ and $_R(ReR)$ are
projective if $\mu$ is bijective. \label{proj-ideal}
\end{Lem}

For a finite-dimensional algebra (over a field $k$) $A$ given by a quiver with relations, any idempotent ideal is generated by an idempotent element. If $e^2=e\in A$ and $eAe\simeq k$, then $Ae\otimes_{eAe}eA\simeq AeA$ implies that $AeA$ is a direct sum of $\dim(eA)$ copies of $Ae$. Note that $AeA$ is the trace of $Ae$ in $_AA$. Thus one can easily decide if $_AAeA$ is projective or not. Let us illustrate this point by an example.
$$A: \quad \xymatrix{
\ar@(ul,dl)[]_{\epsilon}3\bullet & 2\bullet\ar@<2.5pt>[r]^{\beta}\ar[l]_{\delta} &\bullet 1\ar@<2.5pt>[l]^{\gamma}\ar@(ur,dr)[]^{\alpha\; ,}
  & \qquad \gamma\beta=\alpha^2, \; \beta\alpha=\beta\gamma=\delta\epsilon=\epsilon^2=0.}
$$
Here, we write the composition of arrows $\alpha: i\to j$ and $\beta: j\to h$ by $\alpha\beta$: $\alpha$ comes first and then $\beta$ follows.
The indecomposable projective modules can be drawn as follows:
$$
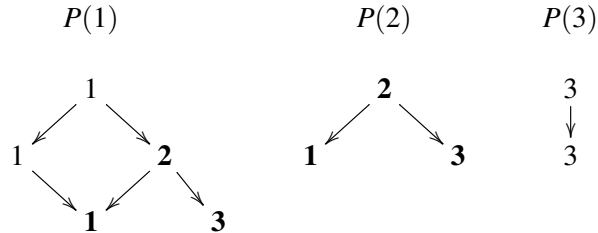
\xymatrix@R=0.4cm@C=0.3cm{
    & P(1)           &          &     &  &     & P(2)            &    & & P(3)         \\
    & 1\ar[dl]\ar[dr]&          &     &  &     & \textbf{2}\ar[dl]\ar[dr]  &   &  &  3\ar[d]    \\
 1 \ar[dr] &                &\textbf{2}\ar[dl] \ar[dr] &     &  & \textbf{1}   &   & \textbf{3}   &                  &3          \\
    & \textbf{1}              &          & \textbf{3}   &  &     &                   &   &  &  \\
}$$
Let $e_i$ be the idempotent of $A$ corresponding to the vertex $i$. Then $_AAe_2A$ is a direct sum of $2$ copies of $Ae_2$ which is denoted by $P(2)$. Note that neither  $_AAe_1A$ nor $_AAe_3A$ are projective since the traces of $Ae_1$ and $Ae_3$ in $_AA=P(1)\oplus P(2)\oplus P(3)$ are not the direct sums of copies of $Ae_1$ and $Ae_3$, respectively.

A variation of Lemma \ref{proj-ideal} is the following result.

\begin{Lem}{\rm \cite{chenxi2016}} Let $R$ be a ring with identity, and let $e$ be an idempotent
element in $R$. Then $_R(ReR)$ (respectively, $(ReR)_R$) is  projective
and finitely generated if and only if $eR(1-e)$ (respectively,
$(1-e)Re$) is  projective and finitely generated as an $eRe$-module
(respectively, a right $eRe$-module), and the multiplication map
${\mu}:(1-e)Re\otimes_{eRe}eR(1-e)\ra (1-e)R(1-e)$ is injective.
 \label{projective}
\end{Lem}

The following result shows that quasi-hereditary rings can be obtained from two smaller quasi-hereditary rings.

\begin{Theo} {\rm \cite{dr1989b}}\label{dr2.4} Let $A$ be a semi-primary ring and $e^2=e\in A$. Then the following are equivalent:

$(1)$ There exists a heredity chain for $A$ containing $AeA$.

$(2)$ Both rings $eAe$ and $A/AeA$ are quasi-hereditary, the multiplication map $$ Ae\otimes_{eAe}eA\lra AeA$$
is bijective, and there exists a heredity chain $\mathscr{I}=(I_i)_{0\le i\le m}$ of $eAe$ such that the
$\mathscr{I}$-filtrations of $(Ae)_{eAe}$ and $_{eAe}(eA)$ are good.

$(3)$ Both rings $eAe$ and $A/AeA$ are quasi-hereditary, the multiplication map $$(1-e)Ae\otimes_{eAe}eA(1-e)\lra (1-e)A(1-e)$$
is bijective, and there exists a heredity chain $\mathscr{I}=(I_i)_{0\le i\le m}$ of $eAe$ such that the
$\mathscr{I}$-filtrations of $((1 - e)Ae)_{eAe}$ and $_{eAe}(eA(1 - e))$ are good.
\end{Theo}

The idempotent ideals in a hereditary chain may not be projective as left or right modules (except the smallest nonzero ideal), but have the following special property.

\begin{Prop}{\rm \cite{dr1989b,P}} If $e$ be an idempotent in a quasi-hereditary ring $A$ such that
$AeA$ belongs to a heredity chain, then the multiplication map $Ae\otimes_{eAe}eA\ra AeA$ is bijective.
\end{Prop}

Note that the lengths of heredity chains of quasi-hereditary algebras may not be maximal with respect to the number of non-isomorphic simple modules. However, we have the following result.

\begin{Prop} {\rm \cite{uy1990}} Any heredity chain of a quasi-hereditary semiprimary ring can be refined to a maximal heredity
chain of the same length as the number of isomorphism classes of simple modules. \label{refine-chain}
\end{Prop}

Finally, we mention that each quasi-heredity ideal $AeA$ of a semiprimary ring $A$ supplies a homological epimorphism from $A$ to $A/AeA$, which provides a recollement of the bounded derived categories of rings. Roughly speaking, a recollement of triangulated categories consists of $3$ triangulated categories and $6$ triangle functors in which there are $4$ adjoint pairs of functors and $3$ fully faithful functors, pointing to the middle triangulated category. For the precise definition of recollements, we refer to \cite{bbd}.

Let $\mathscr{D}^b(A)$ denote the bounded derived category of $A\modcat$. We have the following result.

\begin{Theo}{\rm \cite{ps1988}} Let $e$ be an idempotent in a semiprimary ring $A$ such that $AeA$ is a heredity ideal of $A$. Then there is a recollement of derived categories of rings:

$$\xymatrix@C=1.2cm{\mathscr{D}^b(A/AeA)\ar[r]^-{res.}
&\mathscr{D}^b(A)\ar[r]^-{eA\otimesL-}\ar@/^1.6pc/[l]^-{\rHom(A/AeA,-)}\ar@/_1.6pc/[l]_-{(A/AeA)\otimesL-}
&\mathscr{D}^b(eAe) \ar@/^1.6pc/[l]^-{\rHom(eA,-)}\ar@/_1.6pc/[l]_-{Ae\otimesL-}}$$where $res.$ is the restriction functor induced by the inclusion $\modcat{A/AeA} \hookrightarrow A\modcat$.
\end{Theo}

\medskip
Thus the derived category $\mathscr{D}^b(A)$ of a quasi-hereditary rings can be stratified by a heredity chain. This implies that one may understand the middle category by the information from the outside two categories which are usually ``smaller" than the one in the middle. For example, if $A$ is an Artin algebra, then the number of non-isomorphic simple $A$-modules is the sum of the numbers of non-isomorphic simple modules over the outside two algebras $A/AeA$ and $eAe$; and the Cartan determinant of $A$ is the product of Cartan determinants of $A/AeA$ and $eAe$. For the definition  of Cartan matrices of algebras, we refer to the next Subsection.

\subsection{Examples}
In this subsection, we give examples and non-examples of quasi-hereditary algebras

\begin{Bsp}{\rm The following classes of algebras or rings are quasi-hereditary by \cite[Theorems 1 and 2]{dr1989}.}\end{Bsp}

(1) Let $A$ be a semiprimary ring. Then $A$ is hereditary if and only if
every chain of idempotent ideals of $A$ can be refined to a heredity chain.

(2) A semiprimary ring of global dimension 2 is quasi-hereditary.

A key step to the proof of (2) is based on the following statement: If $e$ is a primitive idempotent of a semiprimary ring $A$ of global dimension at most $2$ such that the Loewy length $L(Ae)$ is
minimal, then $AeA$ is a heredity ideal of $A$.

(3) Let $A$ be a representation-finite Artin algebra. Then the Auslander algebra of $A$ is quasi-hereditary.

This is a consequence of (2), but there was also another proof by splitting filtrations in \cite{dr1989a}.

Recall that the \emph{Auslander algebra} of a representation-finite Artin algebra $A$ is defined as the endomorphism algebra of the direct sum of all non-isomorphic indecomposable $A$-modules. Homologically, an Artin algebra $\Lambda$ is called an \emph{Auslander algebra} if it has global dimension at most $2$ and dominant dimension at least $2$.

\medskip
Given an Artin algebra $A$ with the Jacobson radical $N$ of nilpotency index $n$, consider the $A$-module $M:=\bigoplus_{j=1}^nA/N^j$. Auslander proved that $\End_A(M)$ has global dimension at most $n$ (see \cite{auslander1974}). For such a construction, Dlab and Ringel showed the strong conclusion on $\End_A(M)$.

\begin{Theo} {\rm \cite{dr1989c}} \label{adr-alg} For any Artin algebra $A$ with the Jacobson radical $N$ of nilpotency $n$, the algebra $\End_A(\oplus_{j=1}^n A/N^j)$ is quasi-hereditary.\end{Theo}

This shows that every Artin algebra is of the form $eAe$ with $A$ being quasi-hereditary and $e=e^2\in A$. Equivalently, every Artin algebra can be realized as the endomorphism algebra of a projective module $Ae$ over a quasi-hereditary algebra $A$.

The algebra $\End_A(\oplus_{j=1}^n A/N^j)$ in Theorem \ref{adr-alg} is sometimes called an \emph{Auslander--Dlab--Ringel} algebra in the literature. Some of its properties and structures were studied by Smal$\o$ in \cite{smalo-a, smalo-b}.

Theorem \ref{adr-alg} was further generalized in \cite{linxi}.
A module $M$ is said to be \emph{local} if the sum of all proper submodules of $M$ is a
proper submodule of $M$. So a local module $M$ has a unique maximal submodule which
is the radical of $M$. A module $M$ is said to be \emph{semi-local} if $M$ is isomorphic to the direct sum of finitely many local modules.

\begin{Theo}{\rm \cite{linxi}} Let $A$ be a finite-dimensional algebra with the radical $N$ over an algebraically closed field. If $M$ is a semilocal module with Loewy length $m$, then $\End_A(\bigoplus_{i=l}^m (M/N^iM)$ is quasi-hereditary.
\end{Theo}

\begin{Bsp}{\rm \label{ex2.5}
(1) Schur (or q-Schur) algebras are quasi-hereditary (see \cite{jagreen}).

(2) The following algebra $A$ over a field $k$, given by the quiver and relations, is quasi-hereditary.
$$\xymatrix{\bullet\ar@<-0.4ex>[r]_{\alpha_{1}}_(0){1}_(1){2}
&\bullet\ar@<-0.4ex>[l]_{\beta_{1}}\ar@<-0.4ex>[r]_{\alpha_{2}}_(1){3}
&\bullet\ar@<-0.4ex>[l]_{\beta_{2}}}
\cdots
\xymatrix{\bullet\ar@<-0.4ex>[r]_{\alpha_{n-2}}_(0){}_(1){}
&\bullet\ar@<-0.4ex>[l]_{\beta_{n-2}}\ar@<-0.4ex>[r]_{\alpha_{n-1}}_(1){n}
&\bullet\ar@<-0.4ex>[l]_{\beta_{n-1}}}
$$ $$\beta_{n-1}\alpha_{n-1}=0, \; \alpha_{i}\beta_{i} =\beta_{i-1}\alpha_{\i-1}, 1< i< n.$$

}\end{Bsp}
A heredity chain of ideals of this algebra $A$ is given by
$$0=J_0 \subset Ae_nA\subset A(e_n+e_{n-1})A\subset \cdots \subset A(e_n +\cdots + e_i)A\subset \cdots \subset A(e_n+\cdots + e_n)A=A.$$
This is, in fact, the Auslander algebra of $k[x]/(x^n)$ (see \cite{dr1992}).

Quasi-hereditary algebras appear also in quantum algebras. For example, certain Birman-Murakami-Wenzl algebras (see \cite{xi2000b}).

\medskip
Let us recall the definition of the Cartan matrix of a semiprimary ring $R$. Let $\{e_1, \cdots, e_n\}$ be a complete set of pairwise orthogonal primitive idempotent elements of $R$, and let $S_i:=Re_i/Ne_i$. Then $S_1,\cdots, S_n$ are  a complete list of non-isomorphic simple $R$-modules. Let $c_{ij}$ be the multiplicity  of $S_i$ in a composition series of the indecomposable projective $R$-module $Re_j$, denoted by $[Re_j:S_i]$. The $n$-tuple $(c_{1j},c_{2j},\cdots,c_{nj})$ is called the \emph{dimension vector} of $Re_j$. Then the \emph{Cartan matrix} of $R$ is defined to be $(c_{ij})\in M_n(\mathbb{N})$, denoted by $C(R)$. Thus the $j$-th column of $C(R)$ is the transpose of the dimension vector of $Re_j$.
The determinant of Cartan matrix of a semiprimary ring is called the \emph{Cartan determinant} of the ring.

\begin{Bsp}{\rm \label{ex2.6}
Finite-dimensional cellular algebras over fields of Cartan determinant $1$ are quasi-hereditary (see \cite[Theorem 1.1]{kx1999}). In fact, Example \ref{ex2.5} is a special case of this example.
}\end{Bsp}

For the definition and basic properties of cellular algebras, we refer the reader to \cite{gl1996} and \cite{kx1998}.

As we have seen, algebras of global dimension at most $2$ are always quasi-hereditary. Naturally, one may ask if this is true for algebras of global dimension $3$.  Unfortunately, this is no longer true.
Uematsu and Yamagata gave the first counterexample in \cite{uy1990}. Let $A$ be the algebra given by the following quiver with relations:

\vspace{-1.7cm}
\[A \xymatrix@R=0.5cm@C=0.5cm{
\quad 1\bullet\;\ar@<-0.25em>[dr]_\alpha &  &  & \bullet 5\ar@<0.25em>[dl]^{\delta}\\
& \bullet 3 \ar[r]^\gamma\ar@<-0.25em>[ul]_{\alpha'}\ar@<0.25em>[dl]^{\beta'} & 4\bullet \ar@<0.25em>[ur]^{\delta'}\ar@<-0.25em>[dr]_{\varepsilon'}\\
2\bullet \ar@<0.25em>[ur]^\beta &  & & \bullet6\ar@<-0.25em>[ul]_{\varepsilon}
}\qquad
\begin{matrix}
	\,\\
	\,\\
	\,\\
	\,\\
	\,\\
	\,\\
0=\beta\beta'=\alpha\gamma=\beta\alpha'\\
=\gamma\varepsilon'=\delta\delta'=\varepsilon\delta',\\
\alpha'\alpha=\beta'\beta, \delta'\delta=\varepsilon'\varepsilon.
\end{matrix}
\]
Then $A$ has global  and dominant dimensions $3$, but is not quasi-hereditary.

\subsection{Basic facts on quasi-hereditary algebras}
In this section we summarize some general properties of quasi-hereditary rings.

In the following, $A$ denotes a primary ring, $N$ stands for the Jacobson radical of $A$, and $J$ is an ideal of $A$.  Let $\pd_A(M)$ denote the projective dimension of an $A$-module $M$, and $\gd(A)$ the global dimension of $A$.
The following result is due to Dlab and Ringel, which collects some basic properties of quasi-hereditary algebras or heredity ideals (see \cite{dr1989}).

\begin{Theo}\label{prop-qha}{\rm \cite{dr1989}}
$(1)$ Let $_AJ$ be projective. If $X$ is an $A/J$-module, then
$\pd_A(X)\le 1 + \pd_{A/J}(X).$
Particularly,  $\pd_A(A/J) \le  1.$

$(2)$ Let $J^2=J$ and $_AJ$ be projective.

$\quad$ $(i)$ If $X$ and $Y$ are $A/J$-modules, then
$$ \Ext_A^i(_AX,{}_AY)\simeq \Ext^i_{A/J}(X,Y) \mbox{ for all } i\ge 0.$$
In particular, $\Ext_A^i(A/J,A/J)=0 \mbox{ for } i\ge 1.$

$\quad (ii)$  $\gd(A/J)\le \gd(A)$.

$(3)$ If $JNJ=0$ and $_AJ$ is  projective, then $\gd(A) \le \gd(A/J) + 2.$

$(4)$ If $A$ is a quasi-hereditary algebra with a heredity chain
$$ 0=I_0\subseteq I_1\subseteq \cdots\subseteq I_n=A$$ of length $n$, then $\gd(A)\le 2n-2$ and $L(A)\le 2^n-1$, where $L(A)$ denotes the Loewy length of $A$.
\end{Theo}

Quasi-hereditary algebras have finite global dimension. The upper bound in Theorem \ref{prop-qha}(4) can be generalized in the following way.

\begin{Theo}{\rm \cite[Theorem 5.4]{apt}} Let $e$ be an idempotent element in an Artin algebra $A$. Then $\gd(A)\le \ pd_A(A/AeA)+
\gd(A/AeA) + \gd(eAe)+1.$
\end{Theo}

By the definition of heredity ideals, quasi-heredity of algebras has a close relation with the Jacobson radicals. Moreover, we will see a stronger connection between them (see \cite{xi1993}).

\begin{Prop}
Suppose that  $A$ is a finite-dimensional algebra over an algebraically close field with the Jacobson radical $N$. If there is an integer $n\ge 2$ such that $A/N^n$ is quasi-hereditary, then $A$ itself is quasi-hereditary.
\end{Prop}

Quasi-heredity of algebras is an invariant of Morita equivalences, namely if two semiprimary rings $A$ and $B$ are
Morita equivalent, then $A$ is quasi-hereditary if and only if so is $B$.

This can not be extended to derived equivalences (for definition and related Morita theory of derived categories, see \cite{ happel1988, RickMoritaTh}), as pointed out by Dlab and Ringel in \cite{dr1989}. In fact, they gave a counterexample even by a tilting module (for definition, see Subsection \ref{sect3.2} below). The example is displayed as follows.
$$ A: \quad \xymatrix@R=0.3cm@C=0.4cm{
 & \bullet 1\ar[dr]^-{\beta} &  & \\
3\bullet\ar[ur]^-{\alpha}& & \bullet 2, \ar[ll]^-{\gamma} & \gamma\alpha\beta=\alpha\beta\gamma\alpha=0.\\
}$$
Then the algebra $A$ has global dimension $4$ without any heredity ideals. We denote by $P(i)$ (respectively, $S(i)$) the indecomposable projective (respectively, simple) $A$-module corresponding to the vertex $i$ of the quiver. Let $T:=P(1)\oplus P(2)\oplus S(2)$ and $A_3:=\End_A(T)$. Then $T$ is a tilting $A$-module, and therefore $A$ and $A_3$ are derived equivalent. We represent $A_3$ by the quiver with relations:
$$A_3: \quad \xymatrix{\bullet\ar@<-0.4ex>[r]_{\alpha_{1}}_(0){1}_(1){2}
&\bullet\ar@<-0.4ex>[l]_{\beta_{1}}\ar@<-0.4ex>[r]_{\alpha_{2}}_(1){3}
&\bullet \;,\ar@<-0.4ex>[l]_{\beta_{2}} &
\beta_{2}\alpha_{2}=\alpha_1\alpha_2=\beta_2\beta_1=0, \; \alpha_2\beta_2 =\beta_{1}\alpha_{1}. }$$
 Note that $A_3$ is quasi-hereditary, while $A$ is not quasi-hereditary.

In general, the following algebra is quasi-hereditary and has Loewy length $3$ and global dimension $2n-2$ ($n\ge 1$):
$$A_n: \xymatrix{\bullet\ar@<-0.4ex>[r]_{\alpha_{1}}_(0){1}_(1){2}
&\bullet\ar@<-0.4ex>[l]_{\beta_{1}}\ar@<-0.4ex>[r]_{\alpha_{2}}_(1){3}
&\bullet\ar@<-0.4ex>[l]_{\beta_{2}}}
\cdots
\xymatrix{\bullet\ar@<-0.4ex>[r]_{\alpha_{n-2}}_(0){}_(1){}
&\bullet\ar@<-0.4ex>[l]_{\beta_{n-2}}\ar@<-0.4ex>[r]_{\alpha_{n-1}}_(1){n}
&\bullet\ar@<-0.4ex>[l]_{\beta_{n-1}}}
$$ $$\alpha_{i-1}\alpha_{i} =\beta_{i}\beta_{i-1} =\beta_{n-1}\alpha_{n-1}=0, \; \alpha_{i}\beta_{i} =\beta_{i-1}\alpha_{i-1}, 1< i< n.$$
This algebra is a block of a representation-finite Schur algebra (see \cite{xi1992} or a Temperley-Lieb algebra \cite{westbury, kx1998}).
For a discussion on dominant dimension of (q-)Schur algebras, we refer to the recent article \cite{fhk}.

\medskip
It is known that the Cartan determinant $\det(C(R))$ of the Cartan matrix $C(R)$ of a semiprimary ring $R$ is $\pm 1$ if the global dimension of $R$ is finite. The Cartan determinant conjecture (CDC) excludes the $-1$ case (see \cite{zacharia1983}).

(CDC): If a semiprimary ring $R$ has finite global dimension, then $\det(C(R))=1$.

This conjecture is open up to date. But it holds true for quasi-hereditary rings (see \cite{bf}).

\medskip
If $A$ is a quasi-hereditary algebra over a field, then the dimension of $A$ over the field can be bounded in terms of Cartan invariants of algebras associated with the labelled species of $A$, see \cite{dr1988} for more details.

\subsection{Hochschild homology and cohomology}
Quasi-hereditary algebras have special homological properties (for instance, finite global dimension). In this subsection, we consider their Hochschild homology and cohomology.

Recall that, for an algebra $A$ over a field $k$, the $i$-th \emph{Hochschild homology} $HH_i(A)$ and $i$-th \emph{Hochschild cohomology} $HH^i(A)$ of $A$ can be defined by
$$ HH_i(A):=\Tor_i^{A^{e}}(A,A) \quad \mbox{ and } HH^i(A):=\Ext^i_{A^{e}}(A,A),$$ respectively, where $A^e:= A\otimes_kA^{\opp}$ is the enveloping algebra of $A$.

The low-dimensional Hochschild homology and cohomology have a simple interpretation, namely, $H^0(A)$ is the center of $A$, that is, $H^0(A)=\{a\in A\mid ax=xa, \; \forall \; x\in A\}$; and $H^1(A)$ is isomorphic to Der$(A)/$Inn$(A)$,
where Der$(A)$ is the $k$-space of all $k$-linear derivations of $A$, and Inn$(A)$ stands
for the $k$-space of all inner derivations of $A$, while $H_0(A)$ is the quotient of $A$ modulo the $k$-space $[A,A]$ spanned by all elements of the form $[x,a]:= xa-ax$ for all $x,a\in A$.

\begin{Prop}{\rm  \cite{z}} If $A$ is a quasi-hereditary algebra with $A/N$ separable, then $HH_i(A)=0$ for all $i\ge 1$.
\end{Prop}
But quasi-hereditary algebras may have non-trivial Hochschild cohomologies. For example, for the Temperley-Lieb algebra $A_n$, we have
$$ \dim_kHH^i(A_n)=\begin{cases} n & i=0,\\ 1 & 1\ge i\ge 2n-2,\\ 0 & i\ge 2n-1.\\
\end{cases} $$
(see \cite{penaxi} for calculations).

\subsection{Some generalizations of quasi-hereditary algebras}

Quasi-hereditary algebras were generalized in different directions. The first one is the notions of stratified algebras and standardly stratified algebras (see \cite{wick, cps1996} and \cite{adl1998}). There is also the notion of cellular algebras (see \cite{gl1996}) which capture a class of special quasi-hereditary algebras (see \cite{kx1998,kx1999}). The other direction is given by generalizing the module theoretical description of quasi-hereditary algebras, which we will see in the next section. In this context, one has the notion of standardly stratified systems and homological systems.

In this subsection we recall the definition of standardly stratified algebras (see \cite{wick,cps1996,adl1998}).

\begin{Def} Let $A$ be a finite-dimensional algebra over a field and $J=AeA$ is an idempotent ideal of $A$ generated by an idempotent $e\in A$.

$(1)$ $J$ is called a \emph{stratifying ideal} of $A$ if the following two conditions hold:

$\quad$ $(i)$ The multiplication map $Ae\otimes_{eAe}eA\ra AeA, \; x\otimes y\mapsto xy$, is injective, and

$\quad$ $(ii)$ $\Tor_i^{eAe}(Ae,eA)=0$ for all $i\ge 1$.

$(2)$ The algebra $A$ is called a \emph{stratified algebra} if there is a finite chain:
$$ (*) \quad 0=J_0\subset J_1\subset \cdots\subset J_n=A$$ of idempotent ideals of $A$ such that $J_i/J_{i-1}$ is a stratifying ideal of $A/J_{i-1}$ for all $1\le i\le n$. Such a chain is called a \emph{stratification} of $A$ of length $n$.

$(3)$ The algebra $A$ is called \emph{standardly stratified algebra} if $A$ is stratified algebra with a stratification $(*)$ such that $J_i/J_{i-1}$ is projective over $A/J_{i-1}$ for all $1\le i\le n$.
\end{Def}

Note that $J=AeA$ is  a stratifying ideal of $A$ if and only if $\Ext^i_A(X,Y)\simeq \Ext^i_{A/J}(X,Y)$ for all $X,Y\in A/J\modcat{}$ and all $i\ge 1$.

By definition, quasi-hereditary algebras are standardly stratified algebras, and standardly stratified algebras are stratified algebras. Note that the $\mathbb{R}$-algebra $\mathbb{R}[x]/(x^2)$ is standardly stratified, but not quasi-hereditary. Another example is the algebra $A$ given by the following quiver with a relation:
$$\xymatrix{1\bullet \ar@/_0.5pc/[r]_-{\alpha} &\ar@/_0.5pc/[l]_-{\beta}\bullet 2 ,& \beta\alpha\beta=0. \\}$$
A stratification of $A$ is: $0\subset Ae_2A\subset A$, but any nonzero ideal of $A$ is not heredity. Thus this algebra is standardly stratified, but not quasi-hereditary.
A standardly stratified algebra is quasi-hereditary if it has finite global dimension.

By definition, \emph{any algebra over a field is standardly stratified} with a stratification of length $1$. But this  does not give us any new information on the algebra itself.  More interesting is the case that the length of a stratification is so long as it can be, namely of maximal length.

\medskip
Similar to heredity ideals, we have a recollement of derived module categories for each stratifying ideal $AeA$ of an Artin algebra $A$ (for example, see \cite{cps1996}):

\smallskip
$$\xymatrix@C=0.8cm{\mathscr{D}^b(A/AeA)\ar[r]
&\mathscr{D}^b(A)\ar[r]\ar@/^1.6pc/[l]\ar@/_1.6pc/[l]
&\mathscr{D}^b(eAe). \ar@/^1.6pc/[l]\ar@/_1.6pc/[l]}$$

\bigskip
This implies that the finitistic dimension of a standardly stratified algebra is finite (see also \cite{au}). Here, the \emph{finitistic dimension} is defined to be the supremum of the projective dimensions of finitely generated modules of finite projective dimension. An open conjecture, called the \emph{finitistic dimension conjecture}, states that any Artin algebra has finite finitistic dimension (see \cite{bass} and \cite[Conjecture (11), p.410]{ars}). Note that this conjecture is verified only for a few classes of algebras.

\medskip
There is a large variety of literature on the study of standardly stratified algebras and related algebras. For example, see \cite{adl1998,adl2008,ahlu,cps1996,fm,fkm,mo,pr,xi06}.

\section{Module theoretical approach to quasi-hereditary algebras\label{sect3}}
In this section we survey module theoretical approaches to quasi-hereditary algebras by Ringel and by Dlab and Ringel. The main references will be \cite{ringel-mz}, \cite{dr1992} and \cite{ar1991}. Also, some related topics and results as well as recent developments will be involved.

\subsection{Module theoretical definition of quasi-hereditary algebras}

In this subsection, let $A$ be a finite-dimensional algebra over a field $k$.

Let $\mathcal{C}\subseteq A\modcat$ be a class of $A$-modules. We denote by  $\mathcal{F}(\mathcal{C})$ the category of $A$-modules $M\in A\modcat$ such that $M$ has a finite chain of submodules of $M$:
$$0=M_0\subseteq M_1\subseteq \cdots \subseteq M_t=M$$
with the property that $M_i/M_{i-1}$ is isomorphic to a module in $\mathcal{C}$ for all $1\le i\le t$. Such a chain is called a \emph{$\mathcal{C}$-filtration} of $M$.

Clearly, $\mathcal{F}(\mathcal{C})$ is closed under extensions and direct sums. In general, $\mathcal{F}(\mathcal{C})$ is not closed under direct summands.

Let $\Lambda$ be a partially ordered set such that the non-isomorphic simple $A$-modules are parameterized by $\Lambda$, that is, $\{S(\lambda)\mid \lambda\in \Lambda\}$ is a complete list of non-isomorphic simple $A$-modules. Let $P(\lambda)$ and $I(\lambda)$ denote the projective cover and the injective envelope of $S(\lambda)$, respectively. When we need to indicate which algebra $A$ is concerned, we write $S_A(\lambda)$, $P_A(\lambda)$ and $I_A(\lambda)$ instead of $S(\lambda)$, $P(\lambda)$ and $I(\lambda)$, respectively.

In defining quasi-hereditary algebras in the language of module theory, the so-called standard modules and co-standard modules play an important role.

Let $\Delta(\lambda)$ be the maximal quotient of $P(\lambda)$ with composition factors of the form $S(\mu)$ for $\mu\le \lambda$. Dually, let $\nabla(\lambda)$ be the maximal submodule of $I(\lambda)$, which has composition factors $S(\mu)$ with $\mu\le \lambda$. In the literature, $\Delta(\lambda)$ and $\nabla(\lambda)$ are called \emph{standard}  and \emph{co-standard modules} , respectively, with respect to the partially ordered set $\Lambda$. Let $\Delta:=\{\Delta(\lambda)\mid \lambda\in \Lambda\}$ and $\nabla:=\{\nabla(\lambda)\mid \lambda\in \Lambda\}$. We consider the subcategories $\mathcal{F}(\Delta)$ and $\mathcal{F}(\nabla)$.
Modules in $\mathcal{F}(\Delta)$ are called \emph{$\Delta$-good} (or $\Delta$-filtrated) modules, and modules in $\mathcal{F}(\nabla)$ are called \emph{$\nabla$-good} (or $\nabla$-filtrated) modules.

Roughly speaking, the category $\mathcal{F}(\Delta)$ consists of these $A$-modules that have ``composition factors" $\Delta(\lambda), \lambda\in\Lambda$, and $\mathcal{F}(\nabla)$ consists of these $A$-modules that have ``composition factors" $\nabla(\lambda),\lambda\in\Lambda$. At the moment, whether these ``composition factors" are uniquely determined by the modules themselves is not yet clear.

Let $A^{\opp}$ be the opposite algebra of $A$, and let $D:=\Hom_k(-,k)$ be the usual duality from $\modcat{A}\to A^{\opp}\modcat$.
Then $\{S_{A^{\opp}}(\lambda)=DS_A(\lambda)\mid\lambda\in \Lambda\}$ is a complete list of non-isomorphic simple $A^{\opp}$-modules. Further, $\nabla(\lambda)=D\Delta_{A^{\opp}}(\lambda)$. Thus any statements on standard modules can be interpreted into the ones for co-standards modules.

If $\lambda$ is a maximal element in $\Lambda$, then $\Delta(\lambda)=P(\lambda)$ is projective and $\nabla(\lambda)=I(\lambda)$ is injective. Suppose that $\lambda_1$ is a maximal element in $\Lambda$, we may consider the trace $J_{\lambda_1}$ of $\Delta(\lambda_1)$ in $A$. Clearly, $J_{\lambda_1}$ is an ideal of $A$. Next, we take  a maximal element $\lambda_2$ in $\Lambda\setminus \{\lambda_1\}$, let $J_{\lambda_2}/J_{\lambda_1}$ be the trace ideal of $\Delta(\lambda_2)$ in $A/J_{\lambda_1}$. In this way, we get a chain of ideals of $A$:
 $$ (**)\quad 0\subseteq J_{\lambda_1} \subseteq J_{\lambda_2}\subseteq \cdots\subseteq J_{\lambda_m}=A,$$
such that $J_{\lambda_i}/J_{\lambda_{i-1}}$ is the trace of $\Delta(\lambda_i)$ in the quotient algebra $A/J_{\lambda_{i-1}}$.
This establishes a connection between $\Delta$ and a chain of idempotent ideals of $A$.

\begin{Def}{\rm \cite{dr1992}, \cite{ringel-mz}}\label{def-module}
A finite-dimensional algebra $A$ over a field is said to be \emph{quasi-hereditary} with respect to a partially ordered set $(\Lambda,\le)$ if the following two conditions are satisfied:

$(1)$ $[\Delta(\lambda):S(\lambda)]=1$ for all $\lambda\in \Lambda$.

$(2)$ $P(\lambda)$ has a $\Delta$-filtration (that is, $P(\lambda)\in \mathcal{F}(\Delta))$ for all $\lambda\in \Lambda$.
\end{Def}

In the following, if $A$ is a quasi-hereditary algebra with respect to a partially ordered set $\Lambda$, we often say that $(A,\le)$ is a quasi-hereditary algebra. The condition $[\Delta(\lambda):S(\lambda)]=1$ is equivalent to saying that $\End_A(\Delta(\lambda))$ is a division ring (see \cite[Lemma 1.6]{dr1992}).

By definition, we have the following two basic properties for a quasi-hereditary algebra $A$.

(a) $\Hom_A(\Delta(\lambda),\Delta(\mu))=0$ for $\lambda\not\le \mu$.

(b) $\Ext^1_A(\Delta(\lambda),\Delta(\mu))=0$ if $\lambda\ge \mu$.

(c) $\Hom_A(\Delta(\lambda),\nabla(\mu))\ne 0$ if and only $\lambda=\mu$.

\smallskip
Note that the property (b) shows that the $\Delta$-composition factors of a $\Delta$-good module is uniquely determined. Thus one may speak about the multiplicity of $\Delta(i)$ in a module $M\in \mathcal{F}(\Delta)$, denoted by $[M:\Delta(i)]$. Hence
$\dim_k(M)=\sum_{\lambda\in \Lambda}[M:\Delta(\lambda)]\dim_k\Delta(\lambda)$.

It is known that if $A$ is quasi-hereditary with respect to $(\Lambda,\le)$ then the chain $(**)$ is a heredity chain in the sense of Definition \ref{def2.1}. Conversely, if an algebra $A$ is quasi-hereditary defined by a heredity chain $(*)$, then one may refine the heredity chain to a heredity chain of maximal length by Proposition \ref{refine-chain}. Thus $A$ is quasi-hereditary with the respect to a linear order on the set of representatives of isomorphism classes of simple $A$-modules. From this point of view, the two definitions of quasi-hereditary algebras are equivalent.

\begin{Bsp}{\rm (1) Dlab and Ringel showed that a finite-dimensional hereditary algebra over a field is quasi-hereditary if and only it is quasi-hereditary with respect to any linear order.

(2) Given a finite-dimensional algebra $A$ over a field, there may be different partial orders to make $A$ a quasi-hereditary algebra with different standard modules. Let $A$ be the algebra given by the quiver with relations:
$$
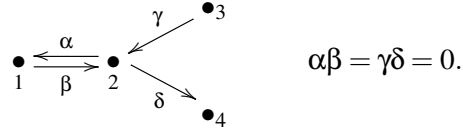
\xymatrix@R=0.2cm{
& & \bullet_{3}\ar[dl]_{\gamma}& \\
\bullet\ar@<-0.4ex>[r]_{\beta}_(0){1}_(1){2} &\bullet\ar@<-0.4ex>[l]_{\alpha}\ar@<0.4ex>[dr]_{\delta}& &\alpha\beta=\gamma\delta=0.\\
 & & \bullet_{4}& \\
}$$
Then the Loewy structure of the indecomposable projective $A$-modules can be pictured visually as follows:
$$
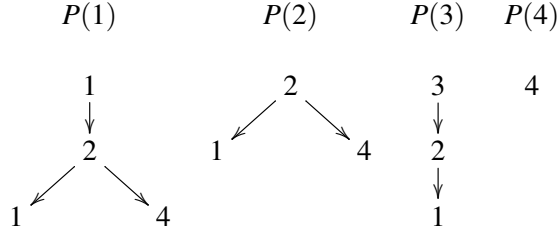
\xymatrix@R=0.4cm@C=0.3cm{
  & P(1)   &          &  &      P(2)        &     & P(3)        & P(4) \\
  & 1\ar[d]&          &  & 2\ar[dl]\ar[dr]  &     &   3\ar[d]   &  4\\
  & 2\ar[dl]\ar[dr]   &  & 1&                  &4    &   2\ar[d]   & \\
 1&        & 4        &  &                  &     &   1         &  \\
}$$

To make $A$ a quasi-hereditary algebra, one may give the partial order $\{ 4>2>1> 3 \}$. Then the  standard modules read as follows:
$$
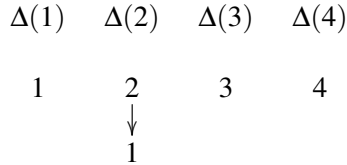
\xymatrix@R=0.4cm@C=0.3cm{
   \Delta(1)      &        \Delta(2) & \Delta(3)        & \Delta(4) \\
     1            &          2\ar[d]  &        3         &  4\\
                  &          1       &                   & \\
}$$
The corresponding heredity chain is of length $4$: $0\subset Ae_4A\subset A(e_4+e_2)A\subset A(e_4+e_2+e_1)A\subset A.$
If one gives the partial order $\{3>4>2>1\}$, then the corresponding standard modules are as follows:
$$
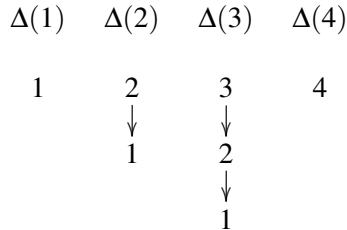
\xymatrix@R=0.4cm@C=0.3cm{
   \Delta(1)      &        \Delta(2) & \Delta(3)        & \Delta(4) \\
     1            &          2\ar[d]  &        3  \ar[d]       &  4\\
                  &          1       &          2 \ar[d]        & \\
                  &                  &         1          & \\
}$$
The corresponding heredity chain is again of length $4$: $0\subset Ae_3A\subset A(e_3+e_4)A\subset A(e_3+e_4+e_2)A\subset A.$ Note that the algebra $A$ has another heredity chain of length $3$: $0\subset A(e_3+e_4)A\subset A(e_3+e_4+e_2)A\subset A.$
 }\end{Bsp}

 We assume that $\Lambda=\{1,2, \cdots, n\}$ with the canonical order. For an $A$-module $M$, we have the dimension vector $\underline{\rm dim}(M)=\big([M:S(1)],[M:S(2)],\cdots, [M:S(n)]\big)\in \mathbb{N}^n$. Let $\underline{\rm dim}(\Delta)$ be the $n\times n$ matrix with $i$-row being $\underline{\rm dim}(\Delta(i))$. Similarly, we have an $n\times n$ matrix $\underline{\rm dim}(\nabla).$ Let $d_\lambda=\dim_k\End_A(\Delta(\lambda))$ for $\lambda\in\Lambda$. The following numerical result holds for quasi-hereditary algebras.

\begin{Prop}{\rm \cite{dr1992}} Let $A$ be a quasi-hereditary algebra over a field $k$ with respect to a partially ordered set $(\Lambda, \le)$. Then

$(1)$ For $\lambda\in \Lambda$ and $M\in\mathcal{F}(\Delta)$, we have $\dim_k\Hom_A(M,\nabla(\lambda))= d_{\lambda}\,[M:\Delta(\lambda)]$.

$(2)$ For $\lambda,\mu\in \Lambda$, we have the Bernstain-Gelfand-Gelfan reciprocity law: $$[P(\mu):\Delta(\lambda)]\, d_{\lambda} =[\nabla(\lambda):S(\mu)]\, d_{\mu}.$$

$(3)$ If $d_\lambda=1$ for all $\lambda\in \Lambda$, then the Cartan matrix $C(A)$ of $A$ has a decomposition: $$C(A)=(\underline{\rm dim}(\Delta))^T\,\underline{\rm dim}(\nabla),$$ where $(\underline{\rm dim}(\Delta))^T$ stands for the transpose of the matrix $\underline{\rm dim}(\Delta)$.
\end{Prop}

Note that (3) follows from (2) by the calculation $$c_{\lambda \mu}=[P(\mu):S(\lambda)]=\sum_{\gamma\in\Lambda}[P(\mu):\Delta(\gamma)][\Delta(\gamma):S(\lambda)]\stackrel{(2)}{=} \sum_{\gamma\in\Lambda}[\nabla(\gamma):S(\mu)][\Delta(\gamma):S(\lambda)].$$

\subsection{$\Delta$-good modules, characteristic tilting modules, and Ringel dualities\label{sect3.2}}
In this section we focus on the categories of good modules over quasi-hereditary algebras.

By specialising results in \cite{ar1991} to quasi-hereditary algebras, Ringel developed a very beautiful theory for the category $\mathcal{F}(\Delta)$ of good modules over a quasi-hereditary algebra in \cite{ringel-mz}.

Throughout this subsection, let $A$ be a finite-dimensional algebra over a field.

\begin{Def} A module $T\in A\modcat$ is called a \emph{tilting $A$-module} {\rm \cite{bb1979, ym86}} if

(i) $\pd_A(T)=n<\infty$,

(ii) $\Ext^i_A(T,T)=0$ for all $i\ge 1$, and

(iii) there is an exact sequence $$0\lra {}_AA\lra T_0\lra\cdots\lra T_n\lra 0$$
with $T_j\in \add(T)$ for all $0\le j\le n$.
\end{Def}

In the case $n=1$, tilting modules are called \emph{classical} tilting modules (see \cite{bb1979}, \cite{hr1982}).

The following result, due to Ringel, provides a strong connection among almost split sequences (see \cite{ars, as, as1981} for definition), tilting modules and quasi-hereditary algebras. This also reflects some natural common features of highest weight categories studied in \cite{cps1988a}.

\begin{Theo}{\rm \cite{ringel-mz}} \label{character} Let $A$ be a quasi-hereditary algebra over a field with respect to a partially ordered set $(\Lambda, \le)$. Then the following hold true.

$(1)$ The $\Delta$-good module category $\mathcal{F}(\Delta)$ has almost split sequences.

$(2)$ For each $\lambda\in \Lambda$, there is an indecomposable $A$-module $T(\lambda)\in \mathcal{F}(\Delta)\cap\mathcal{F}(\nabla)$ such that we have an exact sequence in $\mathcal{F}((\Delta)$:
$$ 0\lra \Delta(\lambda)\lra T(\lambda)\lra X(\lambda)\lra 0,$$
where the $\Delta$-composition factors of $X(\lambda)$ are of the form $\Delta(\mu)$ with $\mu < \lambda$;
and an exact sequence in $\mathcal{F}(\nabla)$:
$$ 0\lra Y(\lambda)\lra T(\lambda)\lra \nabla(\lambda)\lra 0,$$
where the $\nabla$-composition factors of $Y(\lambda)$ are of the form $\Delta(\mu)$ with $\mu > \lambda$.
Moreover, $T:=\bigoplus_{\lambda\in \Lambda} T(\lambda)$ is a tilting $A$-module.

$(3)$ $\mathcal{F}(\Delta)\cap\mathcal{F}(\nabla)=\add(T)$, where $T$ is defined in $(2)$. Moreover, $\mathcal{R}(A):=\End_A(T)$ is a quasi-hereditary algebra with respect to the opposite partially ordered set $\Lambda^{\opp}$ of $\Lambda$. Its standard modules $\Delta_{R(A)}(\lambda)$ are $\Hom_A(T,\nabla(\lambda))$, $\lambda\in \Lambda.$

$(4)$ $\Ext^n_A(\Delta(\lambda),\nabla(\mu))\simeq \begin{cases}k & \mbox{ if } n=0 \mbox{ and  } \lambda=\mu, \\ 0 & otherwise. \end{cases}$
\end{Theo}

\medskip
Following \cite{erdmann1994}, the algebra $\mathcal{R}(A)$ in Theorem \ref{character}(3) is called the \emph{Ringel dual} of the quasi-hereditary algebra $(A,\Lambda)$. Observe that $\mathcal{R}(A)$ has the partially ordered set $\Lambda^{\opp}$ parameterizing its simple modules. If $A$ is a basic quasi-hereditary algebra, then $\mathcal{R}\big(\mathcal{R}(A)\big)\simeq A$.
Thus, for an arbitrary quasi-hereditary algebra $A$, $\mathcal{R}(\mathcal{R}(A))$ is Morita equivalent to $A$.

By a result of Auslander--Reiten in \cite{ar1991}, the characteristic tilting module $T$ determines the categories of $\Delta$-good and $\nabla$-good modules in the follwong way, namely
$$\mathcal{F}(\Delta)=\{M\in\modcat{A} \mid \Ext^i_A(M,T)=0 \; \mbox{ for all }\; i\ge 1\}$$ and
$$\mathcal{F}(\nabla)=\{N\in \modcat{A} \mid \Ext^i_A(T,N)=0 \; \mbox{ for all }\; i\ge 1\}.$$

\medskip
The characteristic tilting modules or their indecomposable direct summands were further extended to algebraic groups in \cite{donkin}. Donkin introduced the notion of tilting modules for rational representations of algebraic groups, where ``tilting modules" are just indecomposable modules which have $\Delta$- and $\nabla$-filtrations. Note that such tilting modules appeared already in the earlier work of Irving \cite{irving}, namely the modules having both a Weyl-module filtration and a Verma-module filtration. For the group $GL_n$, one has to consider the Schur algebras $S(n,r)$. Donkin showed that if $n\ge r$ then Schur algebra  $S(n,r)$ is self Ringel-dual, that is, the Ringel dual of Schur algebra  $S(n,r)$ is isomorphic to itself.

In general, a quasi-hereditary algebra and its Ringel dual may be quite different. In \cite{dengxi96}, it was shown that a Ringel dual may has very large dimension, compared with the given quasi-hereditary algebra. We illustrate this point by an example.

\begin{Bsp}{\rm \label{ex2.10} Dual extension of a quiver with relations. Suppose that $B$ is a finite-dimensional algebra over a field $k$, which is given by a quiver $Q_B=(Q_0, Q_1)$  with relations $\{\rho_i\mid i\in I_B\}$, where $Q_0$ is the vertex set of $Q_B$, $Q_1$ is the arrow set of $Q_B$, and $I_B$ is an index set. We denote by $Q'=(Q_0, Q_1')$ the opposite quiver of $Q_B$. Let $A$ be the algebra given by the quiver $(Q_0,Q_1\cup Q'_1)$ with relations $\{\rho_i\mid i\in I_B\}\cup \{\rho'\mid i\in I_B\}\cup \{\alpha'\beta\mid \alpha, \beta\in Q_1\}$. Then $A$ is a finite-dimensional algebra over $k$, called the \emph{dual extension} of $B$ (see \cite{xi1994}). For more general constructions, we refer to \cite{dengxi1995, xi2000}.

If we take $B$ to be the quiver $$ 1\bullet \longleftarrow \bullet\longleftarrow\cdots\longleftarrow \bullet n, $$then the \emph{dual extension} $A$ of $B$ is quasi-hereditary with the order $1<2<\cdots< n$. Moreover, we have
$$ \dim_k(A)=\frac{1}{6}n(n+1)(2n+1) \mbox{ and } \dim_k(\mathcal{R}(A))= \frac{1}{3}(n^4-1).$$
}\end{Bsp}

For further information on Ringel duality of quasi-hereditary algebras, we refer to the works \cite{bk2018}, \cite{cruz} and\cite{erdmann1994}.
It is shown that the Ringel duality is actually an instance of Koszul duality \cite{bk2018}. For the proof of this result, we refer to the original article \cite{bk2018}.

\subsection{Standardizations}

In this subsection, we present a method, due to Ringel, and to Dlab and Ringel, to produce a quasi-hereditary algebra such that its $\Delta$-good module category is just the prescribed subcategory of modules with certain filtration (see \cite{dr1992}).

Given an abelian $k$-category $\mathcal{A}$ (with $k$ a field),  Dlab and Ringel introduced the notion of standardizable sets in $\mathcal{A}$.

\begin{Def}{\rm \cite{dr1992}} A finite set $\Theta=\{\theta(\lambda)\mid \lambda\in \Lambda\}$ of objects in $\mathcal{A}$ is said to be \emph{standardizable} if the following two conditions are fulfilled:

$(1)$ $\dim_k\Hom_{\mathcal{A}}(\Theta(i),\Theta(j))<\infty$ and $\dim_k\Ext^1_{\mathcal{A}}(\Theta(i),\Theta(j))<\infty$ for all $i,j\in \Lambda$.

$(2)$ The quiver with the vertex set $\Lambda$ and the arrow set $\{i\ra j\mid \rad_{\mathcal{A}}(\Theta(i),\Theta(j))\ne 0 \mbox{ or } \Ext^1_{\mathcal{A}}(\Theta(i),\Theta(j))$ $\ne 0, i,j\in \Lambda\}$, has no oriented cycle.
\end{Def}

Here, by $\rad_{\mathcal{A}}(X,Y)$ we mean the set of morphisms $f:X\to Y$ in $\mathcal{A}$ such that, for any morphism $g: Z\to M$ and $h:N\to Z$ in $\mathcal{A}$, the composition $gfh\in \End_{\mathcal{A}}(Z)$ is not invertible.

The condition (2) in the above definition implies that $\Lambda$ is a partially ordered set: $i\le j$ if there is a finite chain of arrows $i=i_0\to i_1\to\cdots\to i_n=j$.

Standardizable sets can be used to construct quasi-hereditary algebras, as indicated in the following result.

\begin{Theo}{\rm\cite{dr1992}}
If $\Theta$ is a standardizable set in an abelian $k$-category $\mathcal{A}$ with $k$ a field, then there exists a quasi-hereditary $k$-algebra $A$, unique up to Morita equivalence, such that the category $\mathcal{F}(\Theta)$ of $\mathcal{A}$ is equivalent to the category  $\mathcal{F}(\Delta_A)$ of $\Delta_A$-good modules in $A\modcat$.
\end{Theo}

If $\Theta$ is a standardizable set of an abelian $k$-category $\mathcal{A}$, then this set is also a standardizable set of $\mathcal{A}^{\opp}$ with the respect to the opposite partially ordered set. Any subsets of a standardizable set are standardizable.

\subsection{Standard modules of small projective dimension}
Now, we consider some special classes of quasi-hereditary algebras. The first one is quasi-hereditary algebras with all standard modules of smaller projective dimension. In this case, the subcategories $\mathcal{F}(\Delta)$ and $\mathcal{F}(\nabla)$ can be described more explicitly.

Recall that a module is called a \emph{torsionless} module if it is a submodule of a projective module. Dually, a module is called a \emph{divisible} module if it is a homomorphism image of an injective module. If $\mathcal{D}\subseteq \mathcal{C}$ are two additive full subcategories of $A\modcat$, we write $\mathcal{C}/\mathcal{D}$ for the quotient category of $\mathcal{C}$ by $\mathcal{D}$, for which the objects are the same as $\mathcal{C}$, and the Hom-set $\Hom_{\mathcal{C}/\mathcal{D}}(X,Y)$ of objects $X$ and $Y$ is the quotient of $\Hom_A(X,Y)$ modulo all homomorphisms which factorize through a module in $\mathcal{D}$.

For $A$-modules $M$ and $X$, let $tr_M(X)$ denote the \emph{trace} of $M$ in $X$, that is, the sum of the images of all homomorphism from $M$ to $X$. Thus $tr_M(X)$ is a submodule of $X$.

\begin{Theo}{\rm \cite{dr1992}} \label{drms} The following are equivalent for a quasi-hereditary algebra over a field with the standard modules $\Delta$ and co-standard modules $\nabla$.

$(1)$ All standard modules have projective dimension at most $1$.

$(2)$ The characteristic tilting module $T$ has projective dimension at most $1$.

$(3)$ $\mathcal{F}(\nabla)$ is closed under factor modules.

$(4)$ All divisible modules belong to $\mathcal{F}(\nabla)$.
\end{Theo}

The dual statement of Theorem \ref{drms} reads as follows.

\begin{Theo}{\rm \cite{dr1992}} \label{drms-daul} The following are equivalent for a quasi-hereditary algebra over a field with the standard modules $\Delta$ and co-standard modules $\nabla$.

$(1)$ All co-standard modules have injective dimension at most $1$.

$(2)$ The characteristic tilting module $T$ has injective dimension at most $1$.

$(3)$ $\mathcal{F}(\Delta)$ is closed under submodules.

$(4)$ All torsionless modules belong to $\mathcal{F}(\Delta)$.
\end{Theo}

The module $T$ in Lemma \ref{drms}(2) is a classical tilting module in the sense of \cite{bb1979}, \cite{hr1982}). Thus one gets a torsion pair $(\mathcal{G}(T),\mathcal{H}(T))$ in $A\modcat$ from $T$, where
$$\mathcal{G}(T):=\{X\in \modcat{A}\mid \Ext^1_A(T,X)=0\}=\{X\in \modcat{A} \mid Y \mbox{ is generated by } T\},$$ and
$$\mathcal{H}(T)):=\{X\in \modcat{A}\mid \Hom_A(T,X)=0\}=\{X\in \modcat{A} \mid X \mbox{ is cogenerated by } {\rm \tau}(T) \} $$where ${\rm \tau}$ stands for the Auslander--Reiten translation.

Under the condition in Lemma \ref{drms}(1), the equality $\mathcal{F}(\nabla)=\mathcal{G}(T)$ holds. How can we say about $\mathcal{F}(\Delta)$ in this case?

\begin{Theo}{\rm\cite{dr1992}}
Let $A$ be a quasi-hereditary algebra over a field with the standard modules $\Delta=\{\Delta(i)\mid i\in \Lambda\}$.
If all standard modules have projective dimension at most $1$ and all co-standard modules have injective dimension at most $1$, then $\mathcal{F}(\Delta)/\add(T)$ is equivalent to $\mathcal{H}(T)$.
\end{Theo}

Actually, the equivalence functor between $\mathcal{F}(\Delta)/\add(T)$ and $\mathcal{H}(T)$ is given by $tr_T(M)$ which is the trace of $T$ in a module $M$.

\medskip
\begin{Def}{\rm  \cite{ringel2010}} A quasi-hereditary algebra $(A, <)$ is said to be \emph{left strong} if one of the properties in Theorem \ref{drms} holds.
\end{Def}

This class of algebras appeared in Iyama's proof to finiteness of the representation dimensions of Artin algebras, and also in other articles \cite{bhrr,dr1992,hv2003}. Ringel gave a systematic study of such a class of algebras in \cite{ringel2010}. Now, we recall one of his main results on left strong quasi-hereditary algebras.

\medskip
For an Artin algebra $A$ and an $A$-module $X\in A\modcat$, one may consider $X$ naturally as an $A$-$\End_A(X)$-bimodule. Let $\gamma(X)$ be the radical of $\End_A(X)^{\opp}$-module $X$, that is, $\gamma(X)=X\rad(\End_A(X))$. Clearly, $\gamma(X)$ is proper submodule of $_AX$ if $X\ne 0$ by Nakayama Lemma. One may iterate this procedure and define $\gamma^2(X):=\gamma(\gamma(X))$. In general, one defines $\gamma^{\; i+1}(X):=\gamma(\gamma^{\;i}(X))$ for $i\ge 3$.  Then there is a descending chain of submodules:
$$  X=\gamma^0(X)\supset \gamma^1(X) \supset\cdots\supset \gamma\,^t(X)=0.$$
The smallest such  $t$ is denoted by $d(X)$. Clearly, $d(X)\le \ell(X)$, the length of $X$. Warning: $\gamma\;^i(X)\neq X\rad^i(\End_A(X))$. Let $Y:=\bigoplus_{2\le j\le d(X)}\gamma^{\; j}(X)$.

\begin{Theo}{\rm \cite{ringel2010}} For $X\in A\modcat$, the algebra $B:= \End_A(X\oplus Y)$ is left strongly quasi-hereditary and $\gd(B)\le d(X)$.
\end{Theo}

Note that the bound $d(X)$ is much better than the bound given by the length of a heredity chain defining a quasi-hereditary algebra (see Theorem \ref{prop-qha}(4) and \cite{dr2008}).

As pointed out by Ringel in \cite{ringel2010}, the left strong quasi-heredity does not imply the right strong quasi-heredity.

Another special class of quasi-hereditary algebras, called lean quasi-hereditary algebras, were introduced and studied in \cite{adl}. For their precise definition and basic properties, we refer to \cite{adl}.

\subsection{Finiteness of the category of $\Delta$-good modules}
When $\mathcal{F}(\Delta)$ is of finite type (that is, $\mathcal{F}(\Delta)$ has finitely many indecomposable modules, up to isomorphism), is a basic question in the study of quasi-hereditary algebras.

Dlab and Ringel provided a complete answer to this question for the Auslander algebra $\Lambda_n$ of $k[x]/(x^n)$ in \cite{dr1992}.

\begin{Theo}{\rm \cite{dr1992}} Let $k$ be a field and $\Lambda_n$ be the Auslander algebra of $k[x]/(x^n)$. Then $\mathcal{F}(\Delta)$
of $\Delta$-good $\Lambda_n$-modules is of finite type for $n\le 5$ and of tubular type $\tilde{\mathbb{E}}_8$ for $n=6$.
\end{Theo}

In fact, the Auslander--Reiten quiver of $\mathcal{F}(\Delta)$ of $\Lambda_n$ for $n\le 5$ is displayed in \cite{dr1992}. For the definition of tubular types, we refer to \cite{ringel1099}.

The modules in $\mathcal{F}(\Delta)$ of $\Delta$-good $\Lambda_n$-modules  have a quite close relation with the
unipotent radical of a parabolic subgroup $P$ of $GL_m(k)$. It was shown in \cite{hillerohr} that the problem of describing all $\Delta$-filtered $\Lambda_n$-modules is the same as that of describing the conjugation classes of elements in the unipotent radical of a parabolic subgroup $P$ of $GL_m(k)$ under the action of $P$. For $\Delta$-filtered modules without self-extensions, there is a combinatorial description in \cite{bhrr}: Given any $\Delta$-filtered module $X$, there is (up to
isomorphism) a unique $\Delta$-filtered module $Y$ without self-extensions which has the same dimension
vector. Here, an $A$-module $M$ without self-extension means $\Ext^1_A(M,M)=0$.

\medskip
For a hereditary algebra $H$,  Dlab and Ringel showed that for any given order of simple $H$-modules, the algebra $H$ is quasi-hereditary. How can we determine the finiteness of $\mathcal{F}(\Delta)$ of $\Delta$-good $H$-modules? An positive answer can be described by positive  definitiveness of certain quadratic forms. The statement is similar to the representation-finiteness of hereditary algebras. For details, we refer to \cite{deng-xi}.

Suppose that $A$ is a quasi-hereditary algebra such that $\mathcal{F}(\Delta)$ is of finite type, say $X_1, \cdots, X_s$ form a complete list of representatives of isomorphism classes of indecomposable $A$-modules in $\mathcal{F}(\Delta)$. Then one may consider the Auslander algebra $\End_A(\mathcal{F}(\Delta)):=\End_A(\bigoplus_{j=1}^sX_j)$ of $\mathcal{F}(\Delta)$, which is called the $\Delta$-\emph{Auslander algebra} of $\mathcal{F}(\Delta)$.

\begin{Theo}{\rm\cite{xijalgebra}}\label{xi-endo}
If $A$ is a quasi-hereditary algebra with $\mathcal{F}(\Delta)$ of finite type, then the $\Delta$-Auslander algebra $\End_A(\mathcal{F}(\Delta))$ of $\mathcal{F}(\Delta)$ is again quasi-hereditary. Moreover, it is a left QF-3 algebra if and only if the injective direct summand of the characteristic tilting module $T$ in $\mathcal{F}(\Delta)$ cogenerates $T$.
\end{Theo}

Recall that an Artin algebra is called a \emph{left $QT$-3 algebra} if there is a faithful left projective, injective module. This is equivalent to saying that the algebra has a minimal faithful left module (see \cite[p. 40-42]{Tachikawa}).

\subsection{Relations to bocs}

The category $\mathcal{F}(\Delta)$ can be described by direct bocs.  These bocs are
constructed as quotients of dg algebras associated with
the $A_{\infty}$-structure on $\Ext_A^*(\Delta, \Delta)$ (see \cite{kko2014} for more details). The underlying algebra is
an exact Borel subalgebra (see \cite{koenig-mz}). For further investigation of the $A_{\infty}$-structure on the Ext-algebra of standard modules, we refer to \cite{th}. For the theory of bocs, we refer to \cite{cw1988} and \cite{roiter1980}.

Recently, a tame-wild dichotomy theorem for the category $\mathcal{F}(\Delta)$ of $\Delta$-filtered modules is proved in \cite{bps2024}. In fact, the authors of \cite{bps2024} prove the tame-wild dichotomy theorem for a much more general situation of the so-called homological systems including the case of $\mathcal{F}(\Delta)$ as a special example.

\medskip
Motivated by Dlab--Ringel's standardizations, the study of the category $\mathcal{F}(\Delta)$ is generalized to standardly stratified systems in \cite{erd-san2003} and studied in various articles,  and further to homological systems (see \cite{msx} for definition). In \cite{bps2023} the authors associate a differential graded tensor algebra, using the structure of $A_{\infty}$-algebra of
a suitable Yoneda algebra, and use its category of modules to describe the category of filtered modules associated to a
given homological system.

\medskip
Similarly, there is also a module theoretical definition of standardly stratified algebras for a partial order on the set of simple modules. In this case, we have standard modules $\Delta(i)$ and co-standard modules $\nabla(i)$.
The difference is that we remove the restriction: $[\Delta(i):S(i)]=1$ for all $i$ (equivalently, $\End_A(\Delta(i))$ is semisimple for all $i$). Similarly, there is defined the categories $\mathcal{F}(\Delta)$ and $\mathcal{F}(\nabla)$. For further information, we refer to, for example, \cite{bps2023,cps1996,cz2019,fkm,mms,mo,pr,xi06}.
We leave all details and formulations to these original articles.

\section{Decompositions of quasi-hereditary algebras\label{sect3+}}
In this section we introduce the exact Borel subalgebras of quasi-hereditary algebras, due to Steffen Koenig \cite{koenig-mz}. Such a class of quasi-hereditary algebras capture many important algebras like $q$-Schur algebras, Brauer algebras and Temperley-Lieb algebras.

Roughly speaking, quasi-hereditary algebras with exact Borel subalgebras realize their standard modules by induction of simple modules over the subalgebras, or projective modules of directed subalgebras. Recall that a finite-dimensional algebra $B$ is said to be \emph{directed} if the simple $B$-modules (up to isomorphism) can be parameterized by a partially ordered set $(\Lambda,\le)$ such that the indecomposable projective $B$-modules $P(\lambda)$ have the property: $\Hom_B(P(\lambda),P(\mu))=0$ unless $\lambda\le\mu$, and $\End_B(P(\lambda))\simeq \End_B(S(\lambda))$ if $\lambda>\mu$, where $S(\lambda)$ is the simple $B$-module and $P(\lambda)$ is the projective cover of $S(\lambda)$ for $\lambda\in \Lambda$.

\subsection{Exact Borel subalgebras}
We recall the definition of Koenig's exact Borel subalgebras of quasi-hereditary algebras.

\begin{Def}{\rm\cite{koenig-mz}} Let $(A, \le)$ be a quasi-hereditary algebra and $B$ a subalgebra of
$A$. Then $B$ is called an \emph{ exact Borel subalgebra} of $(A,\le)$ if and only if the
following three conditions are satisfied:

$(1)$ $(B,\le)$ is has the same partially ordered set as $A$ for parameterize simple modules, and
$(B, \le)$ is a directed algebra.

$(2)$ $A_B$ is a projective right $B$-module.

$(3)$ $A\otimes_BS(\lambda)=\Delta(\lambda)$ for all $\lambda\in \Lambda$.

\end{Def}

Observe that the exact Borel subalgebra of a quasi-hereditary algebra is the quasi-hereditary analogue of the universal enveloping algebra of a Borel subalgebra of a semisimple complex Lie algebra.

Not every quasi-hereditary algebra has an exact Borel subalgebra. For example, let $A$ be the algebra given by the following quiver with a commutative relation:
$$A: \quad\xymatrix@R=0.2cm@C=0.6cm{
                         &\bullet_{2}\ar[dl]_-{\beta}& &\\
{}_3\bullet&     & \bullet_1\; ,\ar[ul]_-{\alpha}\ar[dl]^{\gamma}&\alpha\beta=\gamma\delta.\\
                         &\bullet_{4}\ar[ul]^-{\delta} & &\\
}$$
For the partial order $\{1<2<3<4\}$, there are $3$ standard modules which are simple corresponding to vertices $1,2,3$, and one standard module is the projective module corresponding to the vertex $4$. It is pointed out in \cite{koenig-mz} that $(A,\le)$ is quasi-hereditary algebra without exact Borel subalgebra.

The existence theorem of Borel subalgebras is given by the following result.

\begin{Theo} {\rm\cite{koenig-mz}} Let $(A, \le)$ be quasi-hereditary with respect to a partially ordered set $\Lambda$, and $B$ a subalgebra of $A$, having
the same simple modules as $A$. Then $B$ is an exact Borel subalgebra of $(A, \le)$ if and only if the following conditions are satisfied:

(1) The algebra $(B,\le)$ is directed.

(2) For each $\lambda\in \Lambda$, there is an isomorphism of $B$-modules by restriction: $\nabla_A(\lambda)\simeq \nabla_B(\lambda)$.
\end{Theo}

Exact Borel subalgebras were discussed for finite-dimensional quasi-hereditary algebras appearing as blocks of category $\mathscr{O}$ of a finite-dimensional semisimple complex Lie
algebra, and for generalized Schur algebras (see \cite{koenig-mz}) and other special classes of algebras \cite{kx1998}. We just mention that Example \ref{ex2.10} provides a class of quasi-hereditary algebras with exact Borel subalgebras: If we start with a directed algebra $B$, then the dual extension of $B$ is quasi-hereditary with $B$ as an exact Borel subalgebra. Note that dual extensions describe the structure of  the Terwilliger algebras of quasi-thin associative schemes \cite{chenxi2025}.

Recently, it is shown that Auslander--Dlab--Ringel algebras over perfect fields have exact Borel subalgebras \cite{cz2019}.

Though quasi-hereditary algebras may have not exact Borel subalgebras in general, we have the following statement when we pass to Morita equivalent algebras.

\begin{Theo}{\rm \cite{kko2014}} Every quasi-hereditary algebra over an algebraically closed field is Morita equivalent to a quasi-hereditary algebra having an exact Borel subalgebra.
\end{Theo}

Thus, when working over algebraically closed fields and considering representations or categorical properties, we can always assume that quasi-hereditary algebras have exact Borel subalgebras.

\subsection{Examples of algebras with a triangular decomposition}
Frequently, quasi-hereditary algebras with exact Borel subalgebras may have triangular decompositions \cite{koenig-mz,koenig-pams}. Recall that an algebra $A$ over a field has a\emph{ triangular decomposition} if there are three subalgebras $S, B,C$ such that

(a) they have the same identity,

(b) $B$ and $C$ are directed algebras such that $B\cap C=S$ is a maximal semisimple subalgebra of $A$,

(c) The multiplication map $B\otimes_SC\to A$ is an isomorphism of $B$-$C$-bimodules.

\medskip
An interesting case of triangular decompositions is $C=B^{\opp}$. There is a construction of such algebras with exact Borel subalgebras in \cite{dyer} and \cite{dengxi1995}. Quasi-hereditary algebras with exact Borel subalgebras often emerge as algebras with a duality fixing simple modules.
In this case, Koenig proved the following result which extends a result on dual extensions in \cite{xi1995}.

\begin{Theo} {\rm\cite{koenig-pams}} Let $(B,\le)$ be a directed algebra. Suppose that $A = B\otimes_S B^{\opp}$ is an algebra
with a triangular decomposition such that $S\simeq B/\rad(B) \simeq A/\rad(A)$. Then $(A,\le )$
is quasi-hereditary. Assume in addition that all standard module $\Delta(\lambda)$ of $A$ are semisimple over $B^{\opp}$. Then
$\gd(A)=2 \cdot \gd(B)$.
\end{Theo}

A special case of triangular decompositions is the following one in which $B$ and $C$ are directed algebras.

\begin{Def}{\rm \cite{xixu}} Let $A$ be an Artin $R$-algebra, $B$, $C$ and $S$ three subalgebras of $A$ $($with
the same identity$)$. We say that $A$ decomposes into a \emph{twisted
tensor product} of $B$ and $C$ over $S$, denoted by $A=B\wedge C$, if

$(1)$ $S$ is a maximal semisimple subalgebra of $A \,(\,$that is,
$S$ is a semisimple $R$-algebra such that $A= S\oplus \rad(A)$ as a
direct sum of $R$-modules $)$ such that $B\cap C=S$.

$(2)$ The multiplication map $\varphi: C\otimes_SB\simeq {}_CA_B$ is
an isomorphism of $C$-$B$-bimodules.

$(3)$ $\rad(B)$\rad$(C)\subseteq$ {\rad }$(C)$\rad $(B),$ where $\rad(B)$ denotes the Jacobson radical of $B$.
\label{decomp}
\end{Def}

For an algebra $A$ that ia a twisted tensor product of $B$ and $C$, we can bound the finitistic dimension of $A$ by the ones of $B$ and $C$ (see \cite[Theorem 1.5]{xixu}).

Examples of algebras that are twisted
tensor products are the dual extensions of directed algebras (see Example \ref{ex2.10}). Other examples may be found in \cite{dengxi1995, kx1998}.

\medskip
Recently, Conde, Dalezios and Koenig investigate triangular decomposition for a new class of quasi-hereditary algebras, called Reedy algebras. They characterise the existence of a Reedy decomposition of an algebra $A$ recursively
via Reedy decompositions of $eAe$ and $A/AeA$, where $e$ is any idempotent of $A$ generating an ideal in
the defining heredity chain (view $A$ as a quasi-hereditary algebra). This is an analogue of Theorem \ref{dr2.4} for Reedy algebras. For more details, we refer to the article \cite{cdk}.

For a discussion of more general cases about exact Borel subalgebras, one is referred to \cite[Theorem 13.10]{bps2023} and \cite{cz2019}. In fact, Bautista, P\'erez and Salmer\'on in \cite{bps2023} extended the existence results of exact Borel subalgebras in \cite{kko2014} to the so-called prestandardly stratified algebras. They showed that every prestandardly stratified algebra over an algebraically closed field is Morita equivalent to a prestandardly stratified algebra having an exact Borel subalgebra. In \cite{cz2019}, the authors extended the notion of exact Borel sublagbras to Borelic pairs for arbitrary finite-dimensional algebras.  They constructed, among others, Borelic pairs for specific families of algebras, including partition algebras, (walled) Brauer algebras, Jones algebras and Temperley-Lieb algebras.

\section{Construction of quasi-hereditary algebras\label{sect4}}
In this section, we mainly explain Dlab--Ringel's construction of quasi-hereditary algebras (or rings).

\subsection{Dlab--Ringel's construction}

The first important construction may be due to Dlab and Ringel in \cite{dr1989b}, which states that all quasi-hereditary algebras with a separable condition can be obtained from a standard procedure.

Given a ring $R$ with identity and a bimodule $_RV_R$, one may form the tensor algebra $\mathcal{T}(_RV_R)$ of $V$ over $R$:
$$ \mathcal{T}(_RV_R)= R + V + V\otimes_RV + V\otimes_RV\otimes V +\cdots V^{\otimes n}+\cdots, $$
where $V^{\otimes n}= \underbrace{V\otimes_RV\otimes_R\cdots\otimes_R V}_n$ is the $n$-fold tensor product of $V$ over $R$. The multiplication on $\mathcal{T}(_RV_R)$ is given in a natural way. Dlab and Ringel constructed a ring $A(\gamma)$ by using this kind of tensor algebras.

Let $C, D$ be rings with identity, $_CS_D, {}_DT_C$ bimodules, and $\gamma: {}_CS_D\otimes_DT_C\to {}_CC_C$ a bimodule homomorphism. Consider the product $C\times D$ of rings and the bimodule $_{C\times D}(S\oplus T)_{C\times D}$ given by
$$ (c,d)(s,t)(c',d'):=(csc',dtd') \mbox{ for } c,c'\in C, d,d'\in D, s\in S, t\in T.$$
Then one can form the tensor algebra $\mathcal{T}(S,T):=\mathcal{T}\big(_{C\times D}(S\oplus T)_{C\times D}\big)$. Let $\mathcal{R}(\gamma)$ be the
ideal of $\mathcal{T}(S, T)$ generated by all elements of the form $s\otimes t - \gamma(s\otimes t)$, $s\in S,t\in T.$
Let $A(\gamma):= \mathcal{T}(S,T)/\mathcal{R}(\gamma)$.

\begin{Theo}{\rm\cite{dr1989b}} \label{dr1989con} Let $C$ and $D$ be quasi-hereditary rings, let $_CS_D, {}_DT_C$ be bimodules, and
$\gamma: {}_CS_D\otimes_DT_C\to {}_CC_C$ be a bimodule homomorphism. If there exists a
heredity chain $\mathscr{J}$ of $C$ such that the $\mathscr{J}$-filtrations of both $_CS$ and $T_C$ are good, then $A(\gamma)$ is a quasi-hereditary ring.
\end{Theo}

Note that the algebra $A(\gamma)$ is Morita equivalent to an algebra $A$ studied in \cite{mv}. Moreover,
Theorem \ref{dr1989con} provides an easy  manner to construct quasi-hereditary rings from given data. A natural and also important question is: Can all quasi-hereditary rings (or algebras) be constructed in this way? The answer is yes. Actually, Dlab and Ringel pointed out that the above procedure can produce almost all quasi-hereditary rings. First, they showed the preliminary proposition.

\begin{Prop}{\rm\cite{dr1989b}}  Let $A$ be a quasi-hereditary ring, and $e\in A$ an idempotent
such that $AeA$ belongs to a hereditary chain of $A$. Assume that there exists
a subring $D$ of $(1 - e)A(1 - e)$ such that $D + (1 - e)AeA(1 - e) =
(1-e)A(1-e)$. Let $C:= AeA, S:= eA(1 - e), T:= (1 - e)Ae$, and $\gamma:
S\otimes_D T\ra C$ be the multiplication map. Then $A = A(\gamma)$.
\end{Prop}

Applying this proposition repeatedly, one gets the following inductive
procedure for constructing quasi-hereditary algebras. For a finite-dimensional algebra $A$, we denote by $s(A)$ the number of isomorphism classes of simple $A$-modules. Recall that, for a quasi-hereditary algebra $A$, we denote by $\mathscr{J}=(J_i)_{0\le j\le m})$ a heredity chain
$$ 0=J_0\subset J_1\subset \cdots\subset J_{m-1}\subset J_m=A.$$

\begin{Theo}{\rm\cite{dr1989b}}
Let $A$ be a nonzero finite-dimensional algebra over a field $k$. Suppose that $A$ is quasi-hereditary with a heredity chain $\mathscr{J} =(J_i)_{0\le j\le m}$. Assume $D:=
A/J_{m-l}$ is a separable algebra over $k$. Then there exists a quasi-hereditary $k$-algebra
$C$ with $s(C)< s(A)$, with a heredity chain $\mathscr{I} = (I_i)_{0\le i\le m-1}$, bimodules $_CS_D,
{}_DT_C$, such that the $\mathscr{I}$-filtrations of $_CS$ and $T_C$ are good, and a bimodule homomorphism $\gamma: {}_CS_D\otimes_DT_C\ra {}_CC_C$ with image contained in $\rad(C)$, such that $A=A(\gamma)$
\end{Theo}

Separability condition can be satisfied for algebras over perfect fields. Thus we have the following result.

\begin{Koro}{\rm\cite{dr1989b}} Let $k$ be a perfect field. Let $A$ be a non-zero quasi-hereditary
finite dimensional $k$-algebra. Then there exists a semisimple $k$-algebra $D$, a
quasi-hereditary $k$-algebra $C$, with $s(C)< s(D)$, and a bimodule homomorphism
$\gamma: {}_CS\otimes_D T_C\ra  {}_CC_C$ such that $A=A(\gamma)$.
\end{Koro}

Thus, by iterating Dlab--Ringel's construction, one can obtain all quasi-hereditary algebras over perfect fields, starting with $C$ the zero algebra and $D$ a semisimple algebra.

\subsection{Dlab--Heath--Marko construction }
Finally, we consider when the endomorphism algebra of a module $M$ over a self-injective algebra $A$ is quasi-hereditary. If the module $M$ contains every indecomposable $A$-module as a direct summand, then $\End_A(M)$ is called a \emph{Morita algebra} in \cite{ky}. This class of algebras is of great interest in understanding homological conjectures in the representation theory and homological algebra of finite-dimensional algebras (see \cite{Tachikawa} and a recent article \cite{cfx}). Now let us mention an interesting construction of quasi-hereditary algebras, due to Dlab, Heath and Marko in \cite{dhm}.

\begin{Theo}{\rm \cite{dhm}}\label{dhm-thm} Let $R$ be a commutative local self-injective $k$-algebra over a splitting
field $k$, $\dim_k(R)= n<\infty$. Let $\mathscr{X}: =\{X(\lambda)\mid \lambda\in \Lambda\}$ be a set of local ideals of $R$ indexed by a
finite partially ordered set $\Lambda$ reflecting inclusions: $X(\lambda')\subset X(\lambda)$ if and only if $\lambda'>\lambda$ in $\Lambda$.
Suppose $R\in \mathscr{X}$. Then $A:=\End_R(\bigoplus_{\lambda\in\Lambda}X(\lambda))$ is a quasi-hereditary algebra
with respect to $\Lambda$ if and only if

$(i)$ $|\mathscr{X}|=n$ and

$(ii)$ $\rad(X(\lambda))=\sum_{\lambda<\mu}X(\mu)$.
\end{Theo}

As pointed out in \cite{dhm}, the algebra $A$ in Theorem \ref{dhm-thm} is an example of lean quasi-hereditary algebras (see \cite{adl} for definition). Moreover, it has a duality fixing simple modules. Other constructions of quasi-hereditary algebras with duality fixing simple modules can be found in \cite{dengxi1995, dyer, xi1994}.

\medskip
{\bf Acknowledgements.} I would like to thank Professor Javad Asadollahi for invitation to write this survey article. The research work was partially supported by the Natural Science Foundation (Grant no. 12226314).

At this opportunity, I would like to heartily thank Claus Michael Ringel who has taught me enjoyable representation theory of algebras, provided me with a lot of chances to learn from him, and built my joyful career life.

\medskip
{\footnotesize Changchang Xi,

School of Mathematical Sciences, Capital Normal University, 100048
Beijing, China, and



{\tt Email: xicc@cnu.edu.cn}
}
\end{document}